\documentclass[12pt]{amsart}

\usepackage{amsmath,amssymb,amsthm,
	amsfonts}

	\newtheorem{dfn}{Definition}[section]
	\newtheorem{thm}[dfn]{Theorem}
	\newtheorem{prop}[dfn]{Proposition}
	\newtheorem{cor}[dfn]{Corollary}
	\newtheorem{lem}[dfn]{Lemma}

	\newtheorem{claim}[dfn]{Claim}

	\newtheorem{case}{Case}

	\numberwithin{equation}{section}

	\newcommand{\dist}{\mathop{\mathit{d}} \nolimits}
	
	\newcommand{\diam}{\mathop{\mathrm{diam}} \nolimits}
	
	\newcommand{\nai}{\mathop{\mathrm{Int}} \nolimits}
	
	\newcommand{\tra}{\mathop{\mathrm{Tra}} \nolimits}

	\newcommand{\lip}{\mathop{\mathcal{L}ip}      \nolimits}
	\newcommand{\me}{\mathop{\mathrm{me}}      \nolimits}
	\newcommand{\supp}{\mathop{\mathrm{Supp}}    \nolimits}

	\newcommand{\obd}{\mathop{\underline{H}_{\lambda}\mathcal{L}\iota_1}    \nolimits}
	\newcommand{\obdd}{\mathop{\underline{H}_{1}\mathcal{L}\iota_1}
	\nolimits}
	\newcommand{\bobd}{\mathop{H_{\lambda}\mathcal{L}\iota_1}
	\nolimits} 
	\newcommand{\bobdd}{\mathop{H_{0}\mathcal{L}\iota_1}        \nolimits}
	\newcommand{\sikaku}{\mathop{\underline{\square}_{\lambda}}                 \nolimits}
	\newcommand{\sikakuu}{\mathop{\underline{\square}_{1}}                 \nolimits}

	\newcommand{\bounasi}{\mathop{\square_{\lambda}}    \nolimits}

\begin{document}

	\title[A note for Gromov's distance functions on the space of mm-spaces]
	{A note for Gromov's distance functions on the space of mm-spaces}
	\author[Kei Funano]{Kei Funano}
	\address{Mathematical Institute, Tohoku University, Sendai 980-8578, JAPAN}
	\email{sa4m23@math.tohoku.ac.jp}
	\keywords{mm-space, box distance function, observable distance function}
	\dedicatory{}
	\date{\today}

	\maketitle

	\setlength{\baselineskip}{5mm}

	\begin{abstract}This is just a note for \cite[Chapter
	 $3\frac{1}{2}_+$]{gromov}. Maybe this note is obvious for a
	 reader who knows metric geometry. I wish that someone study
	 further in this direction.
	\end{abstract}
\hspace{4cm}Comments and questions are welcome.
	\section{The box distance function}
	\begin{dfn}\upshape
	 Let $\lambda \geq 0$ and $(X,\mu)$ be a measure
	   space with $\mu(X)< +\infty$. For two maps $\dist_1,\dist_2:X\times X \to
	   \mathbb{R}$, we define a number $\bounasi (\dist_1,\dist_2)$
	   as the infimum of $\varepsilon >0$ such that there exists a
	   measurable subset $T_{\varepsilon}\subseteq X$ of measure
	   at least $\mu(X)-\lambda \varepsilon$ satisfying
	   $|\dist_1(x,y)-\dist_2(x,y)|\leq \varepsilon$ for any $x,y
	   \in T_{\varepsilon}$.
	 \end{dfn}
	   It is easy to see that this is a distance function on the set of all
	   functions on $X\times X$, and the two distance functions $\square_{\lambda}$ and
	   $\square_{\lambda'}$ are equivalent to each other for any
	   $\lambda,\lambda'>0$. 

An \emph{mm-space} is a triple $(X, \dist_X, \mu_X)$, where $\dist_X$ is a complete separable
	metric on a set $X$ and $\mu_X$ a finite Borel measure on $(X,
	\dist_X)$. Two mm-spaces are \emph{isomorphic} to each other if
there is a measure preserving isometry between the supports of their
measures. We denote by $\mathcal{L}$ the Lebesgue measure on $\mathbb{R}$.
	   \begin{dfn}[parameter]\upshape Let $X$ be an mm-space and
	    $\mu_X(X)=m$. Then, there exists a Borel measurable map
	    $\varphi :[0,m]\to X$ with $\varphi
	    _{\ast}(\mathcal{L})=\mu_X$, where
	    $\varphi_{\ast}(\mathcal{L})$ stands for the push-forward
	    measure of $\mathcal{L}$ by $\varphi$. We call $\varphi$ a
	    \emph{parameter} of $X$. 
	    \end{dfn}
	    Note that if the support of $X$ is not a one-point, then its
	    parameter is not unique.
	   \begin{dfn}[Gromov's box distance function]\upshape If two mm-spaces $X,Y$ satisfy $\mu_X(X)=\mu_Y(Y)=m$, we define
		\begin{align*}
		 \sikaku (X,Y):= \inf \square_{\lambda}(\varphi_X^{\ast}\dist_X,
\varphi_Y^{\ast}\dist_Y),
\end{align*}where the infimum is taken over all parameters 
$\varphi_X:[0,m] \to X, \ \varphi_Y :[0,m] \to Y$, and
$\varphi_X^{\ast}\dist_X$ is defined by $\varphi_X^{\ast}\dist_X(s,t):=\dist_X(\varphi_X(s),\varphi_X(t))$ for $s,t \in [0,m]$.
If $\mu_X(X)<\mu_Y(Y)$, putting $m:=\mu_X(X), m':=\mu_Y(Y)$
, we define
\begin{align*}
\sikaku (X,Y):= \sikaku \Big(  X, \frac{m}{m'}Y \Big)+m'-m, 
\end{align*}where $(m/m')Y:=(Y,\dist_Y,(m/m')\mu_Y)$. 
\end{dfn}

We denote by $\mathcal{X}$ the space of all isomorphic class
of mm-spaces. $\sikaku$ is a distance function on $\mathcal{X}$ for any $\lambda \geq
0$ (See Theorem \ref{kusomoreru}). Note that the distances $\sikaku$ and $\underline{\square}_{\lambda'}$ are
equivalent to each other for distinct $\lambda , \lambda'>0$. 

	The following two lemmas are easy to prove, so we omit the proof.

	\begin{lem}\label{a1}Assume that two mm-spaces $X, Y$ satisfy $m:=\mu_X(X)= \mu_Y(Y)$ and a Borel measurable map $\Phi:[0,m]\to [0,m]$ satisfies
	 $\Phi_{\ast}(\mathcal{L})=\mathcal{L}$. Then, both $\varphi_X \circ {\Phi}:[0,m]\to X$
	 and $\varphi_Y \circ {\Phi}:[0,m]\to Y$ are parameters, and the inequality
	 
\begin{align*}
	  \square_{\lambda}((\varphi_{X}\circ
	  {\Phi})^{\ast}\dist_{X},(\varphi_{Y}\circ
	  {\Phi})^{\ast}\dist_{Y})
	  \leq \square_{\lambda}(\varphi_{X}^{\ast}\dist_{X},\varphi_{Y}^{\ast}\dist_{Y})
	  \end{align*}holds. 
	 \end{lem}

	\begin{lem}\label{abaa}Assume that two mm-spaces $X,Y$ satisfy $m:=\mu_X(X)=\mu_Y(Y)$ and let $0<\alpha \leq 1$. Then, we have
	 \begin{align*}
	  \alpha \sikaku (X,Y)\leq \sikaku (\alpha X, \alpha Y)\leq \sikaku (X,Y). 
	  \end{align*}
	 \end{lem}

The following lemma is the key to prove the triangle inequality for $\sikaku$.
	\begin{lem}\label{oioi}Let $(X,\dist_X,\mu_X)$ be a mm-space and $\varphi_X:[0,m]\to X, \psi_X:[0,m]\to X$ be two parameters. Then, for any $\varepsilon >0$,
	 there exist two Borel measurable maps $\Phi_1, \Phi_2:[0,m]\to [0,m]$ such that $\Phi_{1
	 \ast}(\mathcal{L})=\mathcal{L}, \Phi_{2
	 \ast}(\mathcal{L})=\mathcal{L}$, and
	 \begin{align*}
	  \square_0\big(  (\varphi_X \circ \Phi_1)^{\ast}\dist_X ,(\psi_X \circ \Phi_2)^{\ast}\dist_X   \big)< \varepsilon. 
	  \end{align*}
 \begin{proof}To prove the lemma, we shall approximate $X$ by a
  countable space. For any $\varepsilon>0$, there exists a sequence $\{ X_i \}_{i=1}^{\infty}$ of pairwise disjoint Borel subsets of $X$ such that
	  $X=\bigcup\limits_{i=1}^{\infty} X_i  $ and $\diam X_i <
  \varepsilon$ for each $i\in \mathbb{N}$. Fix a point $x_i\in X_i$
	  for each
	  $i\in \mathbb{N}$. We define
  a distance between $x_i$ and $x_j$ by $\dist_{X'}(x_i,x_j):=\dist_X(x_i,x_j)$, and a Borel measure $\mu_{X'}$ on $X'$ by $\mu_{X'}(\{  x_i  \}):=\mu_{X}(X_i)$.
	   Define two maps $\varphi_{X'}:[0,m)\to X'$ and $\varphi_{X'}:[0,m)\to X'$ by $\varphi_{X'}(t):=x_i$ for $t\in
	  \varphi_X^{-1}(X_i)$ and $\psi_{X'}(t):=x_i$ for $t\in
  \psi_{X}^{-1}(X_i)$. It is easy to see that $\square_0
	  (\varphi_{X}^{\ast}\dist_X, \varphi_{X'}^{\ast}\dist_{X'})< 2\varepsilon$ and
	  $\square_{0}(\psi_{X}^{\ast}\dist_X,\psi_{X'}^{\ast}\dist_{X'})<
  2\varepsilon$. Put $\Phi_{X'}(t):= x_1$ for $t\in \big[0,\mu_{X'}(\{ x_1 \})\big)$, and $\Phi_{X'}(t):=
	  x_i$ for $t\in \Big[\sum\limits_{k=1}^{i-1} \mu_{X'}(\{ x_k
  \}), \sum\limits_{k=1}^{i} \mu_{X'}(\{ x_k \})\Big)$, $i=2,3, \cdots ,n$.

We construct a Borel measurable map $\Phi_1^1
	   :\big[0,\mu_{X'}(\{ x_1 \})\big)\to \varphi_{X'}^{-1}(\{
	   x_1\})$ as follows: There is a sequence $\{  K_n    \}_{n=1}^{\infty}$ of compact subsets of $\varphi_{X'}^{-1}(\{  x_1  \})$ such that
	  $K_1\subseteq K_2 \subseteq \cdots$ and $\mathcal{L} (K_n) \to \mathcal{L}\big( \varphi_{X'}^{-1}(   \{ x_1\})\big)$. Take a Borel measurable map $\Phi_1^{11}:\big[0, \mathcal{L}(K_1)\big)\to K_1$ such that
	  $(\Phi_{1}^{11})_{\ast}(\mathcal{L})=\mathcal{L}$. For each
  $i=2,3,\cdots $, we find a sequence $ \{ (a_n^i,b_n^i)
  \}_{n=1}^{\infty}$ of pairwise disjoint open intervals such that
  $K_i\setminus K_{i-1}= K_i \cap
  \bigcup\limits_{k=1}^{\infty}(a_k^i,b_k^i)$. Take Borel measurable
  maps $\Psi_1:I_1:=\big[\mathcal{L}(K_{i-1}),
	    \mathcal{L}(K_{i-1})+\mathcal{L}(K_i\cap [a_1^1,b_1^i] )\big)\to K_i
	  \cap [a_1^i ,b_1^i]$ and $\Psi_k:I_k:=\big[ \mathcal{L}(K_{i-1})
	    +\sum\limits_{l=1}^{k-1} \mathcal{L}(K_i \cap [a_l^i,b_l^i]),          \mathcal{L}(K_{i-1})
	    +\sum\limits_{l=1}^{k} \mathcal{L}(K_i \cap
	    [a_l^i,b_l^i])\big)\to K_i\cap [a_k^i,b_k^i]$, $k=2,3,
	    \cdots$, such that $(\Psi_k)_{\ast}(\mathcal{L})=\mathcal{L}$
	    for $k=1,2, \cdots$. By modifying each $\Psi_k$, we may assume
	    that $\Psi_k(I_k)\subseteq K_i\cap (a_k^i,b_k^i)$.
	    Then we define a Borel measurable map $\Phi_{1}^{1i}: \big[
	  \mathcal{L}(K_{i-1}),\mathcal{L}(K_i)\big)\to K_i\setminus
	    K_{i-1}$ by $\Phi_{1}^{1i}(t):= \Psi_k(t)$ if $t\in
	    I_k$. Put $\Phi_1^1(t):=\Phi_1^{11}(t)$ for $t\in \big[ 0,
	  \mathcal{L}(K_1)  \big)$ and $\Phi_1^1(t):= \Phi_1^{1i}(t)$ for $t\in \big[  \mathcal{L}(K_{i-1}) , \mathcal{L}(K_i)
	  \big)$. It is obvious that this map $\Phi_1^1$ satisfies
	   $(\Phi_1^1)_{\ast} (\mathcal{L})=\mathcal{L}$. In this way,
	  we find a sequence of Borel measurable maps $\Big\{ \Phi_1^i :   \Big[ \sum\limits_{k=1}^{i-1}\mu_{X'}(   \{   x_k\}),
	  \sum\limits_{k=1}^i \mu_{X'}(\{ x_k \})      \Big) \to
	  \varphi_{X'}^{-1}(\{ x_i  \})
	   \Big\}_{i=2}^{\infty} $ such that $   (\Phi_1^{i})_{\ast} (\mathcal{L})=\mathcal{L}$ for each $i=2,3, \cdots$.

	  Define a Borel measurable map $\Phi_1:[0,m)\to [0,m)$ by $\Phi_{1}(t):= \Phi_1^1(t)$ for $t\in
	  \big[0,\mu_{X'}(\{ x_1\})\big)$ and $\Phi_1(t):= \Phi_1^i(t)$ for $t\in \Big[  \sum\limits_{k=1}^{i-1}\mu_{X'}(\{  x_k
	  \}), \sum\limits_{k=1}^i \mu_{X'}(\{ x_k   \})\Big),  i=2,3, \cdots$. From the above construction, it follows that $\Phi_{1
	  \ast}(\mathcal{L})=\mathcal{L}$ and $\Phi_{X'}=\varphi_{X'}\circ \Phi_1$. In the same way, we find a Borel
	  measurable map $\Phi_2:[0,m)\to [0,m)$ such that $\Phi_{2 \ast}\mathcal{L}=\mathcal{L}$ and $\Phi_{X'}=
	  \psi_{X'}\circ \Phi_2$. Therefore, by using Lemma \ref{a1}, we have
	  \begin{align*}
	   \square_0 \big( (\varphi_{X'}\circ \Phi_1)^{\ast}\dist_{X}, (\psi_X \circ \Phi_2)^{\ast} \dist_X      \big)\leq \ &
	   \square_0\big( (\varphi_X \circ \Phi_1)^{\ast} \dist_X , (\varphi_{X'}\circ \Phi_1)^{\ast} \dist_{X'}       \big)\\
	    & + \square_0 \big( (\psi_{X'}\circ \Phi_2)^{\ast}\dist_{X'}   , (\psi_X \circ \Phi_2)^{\ast} \dist_X \big)\\
	   \leq \ & \square_0(\varphi_{X}^{\ast}\dist_X,\varphi_{X'}^{\ast}\dist_{X'})\\
	   & \hspace{2cm}+\square_0 (\psi_{X'}^{\ast}\dist_{X'},\psi_{X}^{\ast}\dist_{X})\\
	   < \ & 4\varepsilon. 
	   \end{align*}This completes the proof.
	  \end{proof}
	 \end{lem}

	\begin{lem}\label{tri}For any $\lambda\geq 0$, $\sikaku$ satisfies the triangle inequality.
	 \begin{proof}
	  Let $(X,\dist_{X},\mu_{X}),
	  (Y,\dist_{Y},\mu_{Y}),(Z,\dist_{Z},\mu_{Z})$ be mm-spaces and
	  put $m:=\mu_X(X),m':= \mu_Y(Y),m'':= \mu_Z(Z)$.

	  \begin{case}
	   $m=m'=m''$. 
	  \end{case}
	   Let $\varphi_{X}:[0,m]\to X,\ \varphi_{Y}:[0,m]\to Y, \ \psi_{Y}:[0,m]\to
 Y,\ \varphi_{Z}:[0,m]\to Z$ be any parameters. By virtue of Lemma \ref{oioi}, for any $\varepsilon >0$, there exists two Borel measurable maps ${\Phi}_{1}:[0,m]\to [0,m],\
	  {\Phi}_{2}:[0,m]\to [0,m]$
	  such that ${\Phi}_{1\ast}(\mathcal{L})=\mathcal{L},\ {\Phi}_{2\ast}(\mathcal{L})=\mathcal{L}$, and
	  \begin{align*}
	   \square_{\lambda}\big((\varphi_{Y}\circ
	   {\Phi}_{1})^{\ast}\dist_{Y},(\psi_{Y}\circ
	   {\Phi}_{2})^{\ast}\dist_{Y}\big)\leq \square_{0}\big((\varphi_{Y}\circ
	   {\Phi}_{1})^{\ast}\dist_{Y},(\psi_{Y}\circ
	   {\Phi}_{2})^{\ast}\dist_{Y}\big)< \varepsilon.
	   \end{align*}Applying Lemma \ref{a1}, we get
	  \begin{align*}
	   &\square_{\lambda}(\varphi_{X}^{\ast}\dist_{X},\varphi_{Y}^{\ast}\dist_{Y})+ 
	   \square_{\lambda}(\psi_{Y}^{\ast}\dist_{Y},\varphi_{Z}^{\ast}\dist_{Z})\\
	   \geq \ &\square_{\lambda}((\varphi_{X}\circ
	   {\Phi}_{1})^{\ast}\dist_{X},(\varphi_{Y}\circ
	   {\Phi}_{1})^{\ast}\dist_{Y})+\square_{\lambda}((\psi_{Y}\circ
	   {\Phi}_{2})^{\ast}\dist_{Y},(\varphi_{Z}\circ
	   {\Phi}_{2})^{\ast}\dist_{Z})\\
	   \geq \ &\square_{\lambda}((\varphi_{X}\circ
	   {\Phi}_{1})^{\ast}\dist_{X},(\varphi_{Z}\circ
	   {\Phi}_{2})^{\ast}\dist_{Z})-\square_{\lambda}((\varphi_{Y}\circ
	   {\Phi}_{1})^{\ast}\dist_{Y},(\psi_{Y}\circ {\Phi}_{2})^{\ast}\dist_{Y})\\
	   \geq \ &\underline{\square}_{\lambda}(X,Z)-\varepsilon,
	   \end{align*}which shows $\sikaku (X,Y) + \sikaku (Y,Z)\geq
	    \sikaku (X,Z) - \varepsilon$.
	  \begin{case}$m \neq m', m=m'$. 
	   \end{case}If $m<m'$, by Lemma
	  \ref{abaa}, we have
	  \begin{align*}
	   \underline{\square}_{\lambda}(X,Y)+\underline{\square}_{\lambda}(Y,Z)= \
	   &\underline{\square}_{\lambda}\Big(X,\frac{m}{m'}Y\Big)+\underline{\square}_{\lambda}(Y,Z)+m'-m\\
	   \geq \
	   &\underline{\square}_{\lambda}\Big(X,\frac{m}{m'}Y\Big)+\underline{\square}_{\lambda}\Big(\frac{m}{m'}Y,\frac{m}{m'}Z\Big)+m'-m
	   \\
	   \geq \ &\underline{\square}_{\lambda}\Big(X,\frac{m}{m'}Z\Big)+m'-m\\
	   = \ &\underline{\square}_{\lambda}(X,Z).
	   \end{align*}If $m>m'$, we have
	  \begin{align*}
	   \sikaku (X,Z)+\sikaku (Y,Z)= \ & \sikaku \Big( \frac{m'}{m}X,Y    \Big)+\sikaku ( Y,Z )+m-m'\\
	   \geq \ & \sikaku \Big(  \frac{m'}{m}X,Z       \Big)+m-m'\\
	   = \ & \sikaku (X,Z). 
	   \end{align*}

	  \begin{case}
	   $m \neq m', m'\neq m'',m=m''$.
	   \end{case}
	  If $m<m'$, we have
	  \begin{align*}
	   \sikaku (X,Y)+\sikaku (Y,Z)=\sikaku \Big( X,\frac{m}{m'}Y        \Big)+ \sikaku \Big( \frac{m}{m'} Y,Z   \Big)+2(m-m')
	   \geq \sikaku (X,Z).
	   \end{align*}If $m> m'$, applying Lemma \ref{abaa}, we get
	  \begin{align*}
	   \sikaku (X,Y)+\sikaku (Y,Z)=\ &\sikaku \Big(  \frac{m'}{m}X,Y \Big)+ \sikaku \Big( Y,\frac{m'}{m}Z      \Big)+2(m-m')\\
	   \geq \ & \sikaku \Big( \frac{m'}{m}X,\frac{m'}{m}Z          \Big)+ 2(m-m')\\
	   \geq \ & \frac{m'}{m}\sikaku (X,Z)+2(m-m').
	   \end{align*}
	  $m\geq \sikaku (X,Z)$ directly implies that
	  \begin{align*}
	   2(m-m')\geq \Big( 1 - \frac{m'}{m}\Big) \sikaku (X,Z). 
	   \end{align*}Thus, we obtain $\sikaku (X,Y) + \sikaku (Y,Z)\geq \sikaku (X,Z)$.

	  \begin{case}$m \neq m' , m \neq m'', m' \neq m''$.
	   \end{case}If $m<m',m'<m''$, by using Lemma \ref{abaa}, we have
	  \begin{align*}
	   \sikaku (X,Y)+\sikaku (Y,Z) =\ & \sikaku \Big(  X, \frac{m}{m'}Y     \Big) +m'-m + \sikaku \Big(
	   Y,\frac{m'}{m''}Z\Big)+m''-m'\\
	   \geq  \ & \sikaku\Big( X, \frac{m}{m'} Y    \Big) + \sikaku \Big(  \frac{m}{m'}Y ,\frac{m}{m''}Z    \Big)+m''-m\\
	   \geq \ & \sikaku \Big(  X,\frac{m}{m''}Z    \Big)+m''-m\\
	   =\ & \sikaku (X,Z).
	   \end{align*}If $m<m',m''<m', m< m''$, by Lemma \ref{abaa}, we get
	  \begin{align*}
	   \sikaku (X,Y)+\sikaku (Y,Z)=\ &   \sikaku \Big( X, \frac{m}{m'} Y     \Big)+m'-m +\sikaku \Big(  \frac{m''}{m'}Y,Z
	   \Big)+m'-m''\\
	   = \ & \sikaku \Big(  X,\frac{m}{m'} Y         \Big)+\sikaku \Big(  \frac{m''}{m'}Y,Z \Big) + 2m'-m-m''\\
	   \geq \ & \sikaku \Big(X, \frac{m}{m'}Y\Big) + \sikaku \Big(  \frac{m}{m'}Y,\frac{m}{m''}Z     \Big)+m''-m\\
	   \geq \ & \sikaku \Big(  X, \frac{m}{m''}Z \Big)+m''-m\\
	   = \ & \sikaku (X,Z). 
	   \end{align*}We prove the same way for the case of $m< m',
	  m''<m',m''<m$. This completes the proof of Lemma \ref{oioi}. 
	  \end{proof}
	 \end{lem}

		 Let $X$ be a mm-space and $M_r$ be the set of all real $r\times r$ matrices. Then we define a Borel measurable map
	$K_r:X^r \to M_r$ by $K_r(x_1,\cdots ,x_r):= \big(  \dist_X(x_i,x_j)      \big)_{i,j}$, and a Borel measure on $M_r$ by
	$\underline{\mu}_r^{X}:=(K_{r})_{\ast}\big((\mu_X)^r\big)$. 
		\begin{thm}[{mm-Reconstruction theorem, \cite[Section
		 $3\frac{1}{2}.5$, $3\frac{1}{2}.7$]{gromov}}]\label{mm}If two
		 mm-spaces $X, X' $ have $\underline{\mu}_r^
		 {X}=\underline{\mu}_r^{X'}$ for all $r\in \mathbb{N}$,
		 then $X$ and $X'$ are isomorphic to each other. 
		 \end{thm}
A. M. Vershik gave the another proof of the reconstruction thereom in
\cite[Section 2, Theorem]{vershik}. We also refer to \cite[Section 2,
Theorem 2.1]{kondo} for his proof. In \cite{kondo}, T. Kondo generalized
the reconstruction theorem to the space of Borel probability measures
on $\mathcal{X}$. 

		\begin{lem}\label{asedakudesuka}Let $(X,\dist,\mu)$ be a mm-space, and $\varphi_X:[0,m]\to X$ be a parameter of $X$. We set $S:= (
	 [0,m],\varphi_X^{\ast}\dist,\mathcal{L} )$. Then, we have $\underline{\mu}_r^{X}= \underline{\mu}_r^{S}$ for all $r=1,2,\cdots $. 
	 \begin{proof}Let $\varphi:[0,m]^r \to X^r$ be a Borel measurable map defined by $\varphi (t_1,\cdots ,t_r):=\big(
	  \varphi_X(t_1),\cdots,\varphi_X (t_r) \big)$. Obviously, $\varphi_{\ast}(\mathcal{L}^r)=(\mu_X)^r$. Therefore, for
	  any Borel subset $A\subseteq M_r$, we obtain
	  \begin{align*}
	   \underline{\mu}_r^{S}(A)=\ &\mathcal{L}^r(\{(t_1,\cdots ,t_r)\in [0,m]^r \mid ({\varphi
	   }_{X}^{\ast}\dist_X (t_i,t_j))_{i,j} \in A \})\\
	   =\ &\varphi_{\ast}(\mathcal{L}^r)(\{
	   (x_1,\cdots ,x_r)\in X^r \mid (\dist_X (x_i,x_j))_{i,j}\in A \})\\
	   =\ &(\mu_X)^r(\{ (x_1, \cdots ,x_r)\in X^r \mid
	   (\dist_X(x_i,x_j))_{i,j}\in A \})\\
	   =\ &\underline{\mu}_r^X (A). 
	  \end{align*}This completes the proof. 
	  \end{proof}
		\end{lem}

\begin{thm}[Gromov, cf.~{\cite[Section $3\frac{1}{2}.6$ Corollary]{gromov}}]\label{kusomoreru}For any $\lambda\geq 0$, $\sikaku$ is a distance
 fuction on $\mathcal{X}$.
		 \begin{proof}Since $\sikaku$ satisfies the triangle
		  inequality, we only prove that $\sikaku (X,Y)=0$
		  implies $X\cong Y$. Supposing that $\sikaku (X,Y)=0$,
		  we shall show $\underline{\mu}_r^{X} =
		  \underline{\mu}_r^{Y}$ for any $r\in
		  \mathbb{N}$. Then, by Theorem \ref{mm}, we get $X
		  \cong Y$.

		  Since $\sikaku (X,Y)=0$, there exist a sequence $\{
		  \varphi_{X,n} \}_{n=1}^{\infty}$ of parameters of $X$
		  and a sequence $\{ \varphi_{Y,n}   \}_{n=1}^{\infty}$
		  of parameters of $Y$ such that $\bounasi
		  (\varphi_{X,n}^{\ast} \dist_{X}, \varphi_{Y,n}^{\ast}
		  \dist_Y ) \to 0$ as $n\to \infty$. Hence, there exist
		  a sequence $\{\varepsilon_n \}_{n=1}^{\infty}$ of
		  positive numbers and a sequence $\{  Z_n
		  \}_{n=1}^{\infty}$ of Borel subsets of $[0,m]$ such
		  that $\varepsilon_n \to 0$ as $n\to \infty$,
		  $\mathcal{L}(Z_n)\geq m- \lambda \varepsilon_n$, and
		  $|\varphi_{X,n}^{\ast} \dist_{X}(s,t)-
		  \varphi_{Y,n}^{\ast} \dist_{Y} (s,t) |\leq
		  \varepsilon_n$ for any $s,t \in Z_n$. Let $U \subseteq
		  M_r$ be an
		  arbitrary open set and denote by $\dist_{M_r}$ the
		  usual Euclidean distance on $M_r$, that is, 
		  \begin{align*}
		   \dist_{M_r} \big(  (a_{ij})_{i,j}, (b_{ij})_{i,j}
		   \big):= \Big( \sum_{i,j=1}^{r} (a_{ij} -b_{ij})^2 \Big)^{1/2}.
		   \end{align*}
		  Put 
		  \begin{align*}
		   X_{n,\varepsilon}\ &:= \{  (t_1 , \cdots , t_r) \in
		   [0,m]^r \mid \big(\varphi_{X,n}^{\ast}  \dist_{X}
		  (t_i,t_j)     \big)_{i,j} \in U \setminus
		   (M_r\setminus U)_{+\varepsilon}           \},\\
		  Y_n \ &:= \{  (t_1 ,\cdots , t_r) \in [0,m]^r \mid
		   \big( \varphi_{Y,n}^{\ast} \dist_{Y}(t_i,t_j)
		   \big)_{i,j}  \in U       \}. 
		  \end{align*}
		  We take $n_0 \in \mathbb{N}$ such that
		  $\varepsilon_n < \varepsilon /r$ for any $n\geq
		   n_0$. 
		  \begin{claim}\label{kurokki}For any $n\geq n_0$, we
		   have $X_{n, \varepsilon } \subseteq Y_n \cup \big([0,m]^r
		    \setminus (Z_n)^r\big)$.
		    \begin{proof}
		    Take any $(t_1, \cdots , t_r)\in
		    X_{n,\varepsilon}$. If $(t_1, \cdots, t_r)\in
		    (Z_n)^r$, then for any $i,j$ we have
		    \begin{align*}
		     |\varphi_{X,n}^{\ast}\dist_X(t_i,t_j)-\varphi_{Y,n}^{\ast}
		     \dist_Y(t_i,t_j) |\leq \varepsilon_n <  \varepsilon /r,
		     \end{align*}which implies that
		    $\dist_{M_r}\big((\varphi_{X,n}^{\ast}\dist_X(t_i,t_j))_{i,j},
		    (\varphi_{Y,n}^{\ast}\dist_Y (t_i,t_j))_{i,j}\big)<
		    \varepsilon $. Hence, we obtain
		    $(\varphi_{Y,n}^{\ast}\dist_Y(t_i,t_j))_{i,j}\in
		    U$. This completes the proof of the claim.
		    \end{proof}
		   \end{claim}

		  Put
		  $S_n:=([0,m],\varphi_{X,n}^{\ast}\dist_X,\mathcal{L})$
		  and
		  $S_n':=([0,m],\varphi_{Y,n}^{\ast}\dist_Y,\mathcal{L})$
		  and let $m:= \mu_X(X)=\mu_Y(Y)$.
		  Combining Lemma \ref{asedakudesuka} and Claim \ref{kurokki}, for any $n\geq
		  n_0$ we have 
\begin{align*}
\underline{\mu}_r^X(U \setminus (M_r \setminus
 U)_{+\varepsilon})=\underline{\mu}_r^{S_n}(U \setminus (M_r \setminus
 U)_{+\varepsilon})=\mathcal{L}^r(X_{n,\varepsilon}) \leq
\  &\mathcal{L}^r\big( Y_n \cup ([0,m]^r\setminus (Z_n)^r)\big)\\
\leq \ &\mathcal{L}^r(Y_n)+\mathcal{L}^r([0,m]^r \setminus (Z_n)^r)\\
\leq \ &\underline{\mu}_r^{S_n'}(U)+rm^{r-1}\lambda \varepsilon_n \\
= \ &\underline{\mu}_r^{Y}(U)+rm^{r-1}\lambda \varepsilon_n.
\end{align*}
In the above inequality, let first $n \to \infty$ and next $\varepsilon \to
 0$. Then we get $\underline{\mu}_r^{X}(U)\leq
		  \underline{\mu}_r^{Y}(U)$. The same argument shows that
		  $\underline{\mu}_r^{Y}(U)\leq
 \underline{\mu}_r^{X}(U)$, which yields $\underline{\mu}_r^{X}(U)=
 \underline{\mu}_r^{Y}(U)$. This completes the proof of Theorem \ref{kusomoreru}. 
		  \end{proof}
\end{thm}
\section{The observable distance function}

For a measure space $(X,\mu)$ with $\mu(X)< +\infty$, we denote by
$\mathcal{F}(X, \mathbb{R})$ the space of all functions on $X$. Given 
$\lambda \geq 0$ and $f, g \in \mathcal{F}(X, \mathbb{R})$, we put 
\begin{align*}
 \me_{\lambda}(f,g):= \inf \{  \varepsilon >0 \mid \mu 
 \big(\{ x\in X \mid |f(x)-g(x)|\geq \varepsilon \} \big)\leq \lambda \varepsilon             \}.
\end{align*}Note that this $\me_{\lambda}$ is a distance function on $\mathcal{F}(X,
 \mathbb{R})$ for any $\lambda \geq 0$ and its topology on
 $\mathcal{F}(X, \mathbb{R})$ coincides with the topology of the 
convergence in measure for any $\lambda >0$. Also, the distance functions
 $\me_{\lambda}$ for all $\lambda >0$ are mutually equivalent. 

 We recall that the \emph{Hausdorff distance} between two closed
subsets $A$ and $B$ in a metric space $X$ is defined by
 \begin{align*}
\dist_H (A,B):= \inf \{  \varepsilon >0 \mid A \subseteq B_{\varepsilon} ,
  B \subseteq A_{\varepsilon}       \},
 \end{align*}where $A_{\varepsilon}$ is a closed $\varepsilon$-neighborhood of $A$.

Let $(X,\mu)$ be a measure space with $\mu(X)< +\infty$. For a
semi-distance function $\dist$ on $X$, we indicate by $\lip_1 (\dist)$ the space of
all $1$-Lipschitz functions on $X$ with respect to $\dist$. Note that
$\lip_1 (\dist)$ is a closed subset in
$(\mathcal{F}(X,\mathbb{R}),\me_{\lambda})$ for any $\lambda\geq 0$.

 \begin{dfn}\upshape For $\lambda \geq 0$ and two semi-distance functions $\dist , \dist'$ on $X$, we define
  \begin{align*}
   \bobd (\dist ,\dist'):= \dist_{H} \big(  \lip_1(\dist ), \lip_1 (\dist')      \big),
   \end{align*}where $\dist_H$ stands for the Hausdorff distance
  function in $(\mathcal{F}(X,\mathbb{R}),\me_{\lambda}).$
  \end{dfn}
This $\bobd$ is actually a distance function on the space of all
semi-distance functions on
$X$ for all $\lambda \geq 0$, and the two distance functions $\bobd$ and
$H_{\lambda'}\mathcal{L}\iota_{1}$ are equivalent to each other for any
$\lambda, \lambda' >0$.
  \begin{lem}\label{mikitii}For any two semi-distance functions $\dist , \dist'$ on $X$, we have
   \begin{align*}
    \bobd (\dist , \dist')\leq \bounasi (\dist, \dist').
    \end{align*}
   \begin{proof}For any $\varepsilon>0$ with $\bounasi (X,Y)<
    \varepsilon$, there exists a measurable subset $T_{\varepsilon} \subseteq
    X$ such that $\mu (X\setminus T_{\varepsilon})\leq \lambda
    \varepsilon$ and $|\dist(x,y)-\dist'(x,y)|\leq \varepsilon$ for any
    $x,y \in T_{\varepsilon}$. Given arbitrary $f\in \lip_1(\dist)$, we
    define $\widetilde{f}\in \mathcal{F}(X, \mathbb{R})$ by $\widetilde{f}(x):= \inf
    \{ f(y)+\dist' (x,y) \mid y\in T_{\varepsilon}  \}$. We see
    easily that $\widetilde{f} \in \lip_1 (\dist')$ and $\widetilde{f}(x)
    \leq f(x)$ for any $x\in T_{\varepsilon}$. Taking any $x\in
    T_{\varepsilon}$, we have
    \begin{align*}
     |f(x)- \widetilde{f}(x)| = \ &f(x)-\widetilde{f}(x)\\ = \ &  \sup
     \{ f(x)-f(y)-d'(x,y) \mid 
     y\in T_{\varepsilon} \}\\  \leq \ &\sup \{
     d(x,y)-d'(x,y) \mid y\in T_{\varepsilon}  \} \\ \leq  \
     &\varepsilon.  
     \end{align*}Therefore, we get $\me_{\lambda}(f, \widetilde{f})\leq
    \varepsilon $, which implies $\lip_1 (\dist) \subseteq \big( \lip_1
    (\dist')      \big)_{\varepsilon}$. Similary, we also have
    $\lip_1(\dist') \subseteq \big( \lip_1 (\dist)
    \big)_{\varepsilon}$, which yields $\bobd (\dist , \dist')\leq
    \varepsilon$. This completes the proof.
    \end{proof}
   \end{lem}

  \begin{dfn}[Observable distance function]If two mm-spaces $X,Y$ satisfy $\mu_X
   (X) =\mu_Y(Y)=m$, we define
   \begin{align*}
    \obd (X,Y):= \inf \bobd (\varphi_{X}^{\ast}\dist_X, \varphi_{Y}^{\ast}\dist_Y),
   \end{align*}where the infimum is taken over all parameters
   $\varphi_X :[0,m]\to X, \ \varphi_{Y}:[0.m]\to Y$. If $\mu_X (X) <
   \mu_Y (Y)$, putting $m:= \mu_X (X), m':= \mu_Y (Y)$, we define 
   \begin{align*}
    \obd (X,Y):= \obd \Big( X, \frac{m}{m'}Y \Big) + m'-m.
    \end{align*}
  \end{dfn}
$\obd$ is a distance function on $\mathcal{X}$ for any $\lambda\geq 0$
(See Theorem \ref{darlin}).
Note that the distance functions $\obd$ and
$\underline{H}_{\lambda'}\mathcal{L}\iota_1$ are equivalent to each
other for any $\lambda, \lambda' >0$.

The proofs of following four lemmas are easy.
\begin{lem}\label{amenotihare}For any parameter $\varphi_X :[0,m]\to X$ of $X$, we have
 \begin{align*}
  \lip_1 (\varphi_X^{\ast}\dist_X)=\{ f\circ \varphi_X \mid f \in \lip_1(\dist_X)    \}.
  \end{align*}
\end{lem}

 \begin{lem}\label{umi1}Assume that two mm-spaces $X,Y$ satisfy $m:=
   \mu_X(X)=\mu_Y(Y)$ and a Borel measurable map
   ${\Phi}:[0,m]\to [0,m]$ satisfies ${\Phi}_{\ast}
   (\mathcal{L}) =\mathcal{L}$. Then, we have
   \begin{align*}
    \bobd \big( (\varphi_X \circ {\Phi})^{\ast} \dist_X,
    (\varphi_Y \circ {\Phi})^{\ast}  \dist_Y
    \big) 
    = \bobd (\varphi_{X}^{\ast} \dist_X , \varphi_{Y}^{\ast} \dist_Y).
    \end{align*}
   \end{lem}

   \begin{lem}\label{umi2}Assume that two mm-spaces $X,Y$ satisfy $m:= \mu_X (X) =
    \mu_Y (Y)$ and let $0 < \alpha \leq 1$. Then, we have
    \begin{align*}
     \alpha \obd (X, Y)\leq \obd (\alpha X, \alpha Y)\leq \obd (X, Y).
     \end{align*}
    \end{lem}

  \begin{lem}\label{umi3}Let X be a mm-space and $\varphi_X :[0,m] \to X,  \psi_X :
   [0,m]\to X$ be two parameters. Then, for any $\varepsilon >0$, there
   exist two Borel measurable maps ${\Phi}_1,
   {\Phi}_2:[0,m] \to [0,m]$ such that $({\Phi}_{1})_{\ast}(\mathcal{L}) = \mathcal{L}$, $({\Phi}_{2})_{\ast}
   (\mathcal{L} )=\mathcal{L}$, and 
   \begin{align*}
    \bobdd \big( (\varphi_X \circ {\Phi}_1)^{\ast}
    \dist_X  ,(\psi_X \circ {\Phi}_2)^{\ast} \dist_X
    \big)< \varepsilon.
    \end{align*}
    \end{lem}

    \begin{thm}[{Gromov, cf.~\cite[Section $3\frac{1}{2}.45$]{gromov}}]\label{darlin}For any $\lambda\geq 0$, $\obd$ is a distance function on
     $\mathcal{X}$.
     \begin{proof}Combining Lemma \ref{umi1}, \ref{umi2}, and
      \ref{umi3}, we see that $\obd$ satisfies the triangle inequality
      in the same way of the proof of Lemma \ref{tri}.

      To prove ``$\obd (X, Y)=0 \Rightarrow X \cong Y$'', we shall approximate
      each $X$ and $Y$ by finite spaces. Take an arbitrary $\varepsilon
      >0$. Then, there exists sequences $\{ X_i \}_{i=1}^{\infty}$, $\{
        Y_j  \}_{j=1}^{\infty}$ of pairwise disjoint Borel subsets of $X$, $Y$ such that
      \begin{itemize}
	    \item[$(1)$]$X=\bigcup\limits_{i=1}^{\infty}X_i$ and $\diam X_i
			\leq \varepsilon$ for any $i \in \mathbb{N}$,
	    \item[$(2)$]$Y= \bigcup\limits_{j=1}^{\infty}Y_j $ and $\diam Y_j
			\leq \varepsilon$ for any $j\in \mathbb{N}$.
      \end{itemize}
      Put $m:= \mu_X(X)=\mu_Y (Y)$. Then, there exists $m_0 \in \mathbb{N}$ such that
      \begin{align*}
       m-\varepsilon \leq \mu_X \Big(  \bigcup_{i=1}^{m_0}X_i
       \Big), m-\varepsilon \leq \mu_Y \Big(  \bigcup_{j=1}^{m_0} Y_j                \Big).
       \end{align*}Since $\obdd (X,Y)=0$, there exist a sequence $\{
      \varepsilon_n \}$ of positive numbers and sequences $\{
      \varphi_{X,n}  \}_{n=1}^{\infty}$, $\{ \varphi_{Y,n}
      \}_{n=1}^{\infty} $ of parameters of $X$, $Y$ such that
      $H_1\mathcal{L}{\iota_1}(\varphi_{X,n}^{\ast}\dist_X,
      \varphi_{Y,n}^{\ast}\dist_Y)< \varepsilon_n$ and $\varepsilon_n
      \to 0$ as $n\to \infty$. For each $i,j=1,\cdots ,m_0$, we fix 
      points $x_i \in X_i$ and $y_j \in Y_j$. Define a function
      $g_{ni}:[0,m]\to \mathbb{R}$ by
      $g_{ni}(s):=\dist_X(\varphi_{X,n}(s),x_i)$ for each $i=1,2,
      \cdots, m_0$. From Lemma \ref{amenotihare}, we have 
      $g_{ni}\in \lip_1(\varphi_{X,n}^{\ast}\dist_X)$. Hence, there
      exists $h_{ni}\in \lip_1(\dist_Y)$ such that $\me_1(g_{ni},h_{ni}\circ
      \varphi_{Y,n})< \varepsilon_n$. Putting 
      \begin{align*}
       A_{ni}:= \{ s\in [0,m] \mid |g_{ni}(s)-(h_{ni}\circ
       \varphi_{Y,n})(s)|< \varepsilon_n            \},
       \end{align*}we get $\mathcal{L}(A_{ni})\geq m-\varepsilon_n$. For
      each $j=1,2, \cdots ,m_0$, we define a function
      $\widetilde{h}_{nj}:[0,m]\to \mathbb{R}$ by
      $\widetilde{h}_{nj}(s):= \dist_Y(\varphi_{Y,n}(s),y_j)$. By the same
      argument as above, there exists $\widetilde{g}_{nj}\in
      \lip_1(\dist_X)$ such that $\mathcal{L}(B_{nj})\geq
      m-\varepsilon_n$, where
      \begin{align*}
       B_{nj}:= \{ s\in [0,m] \mid
 |\widetilde{h}_{nj}(s)-(\widetilde{g}_{nj}\circ \varphi_{X,n})(s)|<
 \varepsilon_n\}.
       \end{align*}
      So, putting
      \begin{align*}
       Z_n:=\varphi_{X,n}^{-1}\Big(\bigcup_{i=1}^{m_0}X_i\Big)
       \cap \varphi_{Y,n}^{-1}\Big( \bigcup_{j=1}^{m_0}Y_j \Big)\cap
       \bigcap_{k=1}^{m_0}A_{nk} \cap \bigcap_{l=1}^{m_0}B_{nl},
      \end{align*}we obtain $\mathcal{L}(Z_n)\geq  2\varepsilon +2m_0
      \varepsilon_n$.

      For any $s,t \in Z_n$, there exist $1\leq i_1,j_1,i_2,j_2\leq m_0$ such
      that 
      \begin{align*}
&s\in \varphi_{X,n}^{-1}(X_{i_1})\cap\varphi_{Y,n}^{-1}(Y_{j_1})\cap
 \bigcap_{k=1}^{m_0}A_{nk}\cap \bigcap_{l=1}^{m_0}B_{nl} \\
\text{and } \ &t\in \varphi_{X,n}^{-1}(X_{i_2})\cap\varphi_{Y,n}^{-1}(Y_{j_2})\cap
 \bigcap_{k=1}^{m_0}A_{nk}\cap \bigcap_{l=1}^{m_0}B_{nl}.
       \end{align*}Since $t\in \varphi_{X,n}^{-1}(X_{i_2})$ and $\diam
      X_{i_2} \leq \varepsilon$, we have 
      \begin{align*}
\dist_X(\varphi_{X,n}(s),\varphi_{X,n}(t))\leq \
 &\dist_X(\varphi_{X,n}(s),x_{i_2})+\dist_X(x_{i_2},\varphi_{X,n}(t))\\
\leq \ &\dist_X(\varphi_{X,n}(s),x_{i_2})+\varepsilon.
       \end{align*}We also get $\dist_X(\varphi_{X,n}(s),x_{i_2})\leq (h_{ni_2}\circ
 \varphi_{Y,n})(s)+ \varepsilon_n$ by $s\in
      \bigcap\limits_{k=1}^{m_0}A_{nk}\subseteq A_{ni_2} $. Therefore,
      we obtain
      \begin{align*}
       \dist_X(\varphi_{X,n}(s),\varphi_{X,n}(t))\leq \ &(h_{ni_2}\circ
       \varphi_{Y,n})(s)+\varepsilon_n +\varepsilon\\
       \leq \ &|(h_{ni_2}\circ \varphi_{Y,n})(s)-(h_{ni_2}\circ
       \varphi_{Y,n})(t)|+|(h_{ni_2}\circ \varphi_{Y,n})(t)|\\
       &  \hspace{8cm}+\varepsilon_n + \varepsilon\\
       \leq \ &
       \dist_Y(\varphi_{Y,n}(s),\varphi_{Y,n}(t))+|(h_{ni_2}\circ 
       \varphi_{Y,n})(t)|+\varepsilon_n+\varepsilon.
       \end{align*}Since $t\in \bigcap\limits_{k=1}^{m_0}A_{nk}\cap
      \varphi_{X,n}^{-1}(X_{i_2})$ and $\diam X_{i_2}\leq \varepsilon$,
      we have $g_{ni_2(t)}\leq \varepsilon$ and $|g_{ni_2}(t)-(h_{ni_2}
      \circ \varphi_{Y,n})(t)|< \varepsilon_n$, and thus $|(h_{n
      i_2}\circ \varphi_{Y,n})(t)|< \varepsilon_n +
      \varepsilon$. Therefore, we obtain
      \begin{align*}
       \dist_X(\varphi_{X,n}(s),\varphi_{X,n}(t))\leq 
      \dist_Y(\varphi_{Y,n}(s), \varphi_{Y,n}(t))+2\varepsilon_n 
      +2\varepsilon.
       \end{align*}A similar argument shows that
      \begin{align*}
       \dist_Y(\varphi_{Y,n}(s),\varphi_{Y,n}(t))\leq
       \dist_X(\varphi_{X,n}(s),\varphi_{X,n}(t))+2\varepsilon_n
       +2\varepsilon.
       \end{align*}Hence, we get
      \begin{align*}
       |\dist_X(\varphi_{X,n}(s),\varphi_{X,n}(t))-\dist_Y(\varphi_{Y,n}(s),\varphi_{Y.n}(t))|\leq
 2\varepsilon_n +2\varepsilon.
       \end{align*}Therefore, we obtain
      \begin{align*}
       \underline{\square}_1(X,Y)\leq
 {\square}_1(\varphi_{X,n}^{\ast}\dist_X,\varphi_{Y,n}^{\ast}\dist_Y)\leq
 2\varepsilon + 2m_0\varepsilon_n.
       \end{align*}So, we get $\sikakuu (X,Y)=0$ and $X\cong Y$. This
      completes the proof.
      \end{proof}
     \end{thm}
     Modifying the proof of Theorem \ref{darlin}, we get the following
     corollary:

     \begin{cor}For any two mm-spaces $X$ and $Y$, we have
      \begin{align*}
       \underline{H}_0 \mathcal{L}{\iota}_1 (X,Y)\leq
       \underline{\square}_0(X,Y) \leq 2\underline{H}_0
       \mathcal{L}{\iota}_1 (X,Y).
       \end{align*}
      \end{cor}

We also refer to \cite[Section $7.4$]{pestov}.
 \section{Another natural method}   

Let $\lambda \geq 0$ and $\varepsilon>0$. A map from an mm-space to a metric space, say
$f:X\to Y$ is called \emph{$\lambda$-Lipschitz up to $\varepsilon$} if 
\begin{align*}
\dist_Y \big(f(x),f(x')\big) \leq \lambda \dist_X (x,x') +\varepsilon
\end{align*}for all $x,x'$ in a Borel subset $X_0\subseteq X$ with
$\mu_X(X\setminus X_0)\leq \varepsilon$. 

\begin{prop}[{cf.~\cite[Section $3\frac{1}{2}.15$, $(3_b)$]{gromov}}]\label{kumogaippai}Let $(X,\dist_X, \mu_X)$, $(Y,\dist_Y, \mu_Y)$ be mm-spaces
 and $\lambda\geq 0$. Let $\varepsilon_n>0$ and $f_n:X \to Y$
 a $\lambda$-Lipschitz up to $\varepsilon_n$ Borel merasurable map and
 assume that $\varepsilon_n \to 0$ as $n\to \infty$ and the sequence $\{
 (f_n)_{\ast}(\mu_X) \}_{n=1}^{\infty}$ converges weakly to
 $\mu_Y$. Then, the sequence $\{ f_n \}_{n=1}^{\infty}$ has a
 $\me_1$-convergent subsequence.
\begin{proof}Without loss of generality, we may assume that $X=\supp
 \mu_X$ and $\mu_X (X)=\mu_Y(Y)=1$. 

By choosing a subsequence, we have
 $\sum\limits_{n=1}^{\infty} \varepsilon_n < +\infty$. From the
 assumption, there exists a Borel subset $X_n\subseteq X$ such that
 $\mu_X(X\setminus X_n)\leq \varepsilon_n$ and $\dist_Y \big( f_n(x),
 f_n(y)      \big)\leq \lambda \dist_X (x,y)+\varepsilon_n$ for any
 $x,y\in X_n$. Put $X_0:= \bigcup\limits_{n=1}^{\infty}
 \bigcap\limits_{i=n}^{\infty} X_i$. Since 
\begin{align*}
\mu_X (X\setminus X_0) \leq \sum_{i=n}^{\infty}\mu_X(X\setminus X_i)\leq
 \sum_{i=n}^{\infty}\varepsilon_i \to 0 \ \text{as } n\to \infty, 
\end{align*}we have $\mu_X(X_0)=1$. Take a countable dense subset $\{p_j
 \}_{j=1}^{\infty} \subseteq X_0$.
 \begin{claim}\label{mikihayoi}The sequence $\{f_{n}(p_1)\}_{n=1}^{\infty}$ has a
  convergent subsequence.
  \begin{proof}The proof is by contradiction. If the sequence $\{
   f_{n}(p_1)\}_{n=1}^{\infty}$ has no convergent subsequence, then the
   subset $A:=\{ f_{1}(p_1), f_2(p_1), \cdots \}$ is a closed subset in
   $Y$, especially, $A$ is
   complete. From the assumption, this set $A$ is not
   compact. Hence, $A$ is not totally bounded, that is, there exists
   $\delta>0$ such that $A$ has no finite 2$\delta$-net. Therefore, by
   choosing a subsequence, we get $B_Y(f_j(p_1),\delta)\cap
   B_Y(f_k(p_1), \delta)=\emptyset$ for any $ j,k$ with $j\neq k$. Take
   $\delta'>0 $ such that $0< \delta' < \delta$ and $\mu_Y\big( \partial B_Y (f_j(p_1),\delta')
   \big)=0$ for any $j\in \mathbb{N}$. Since
   \begin{align*}
    \mu_Y \Big( \partial \Big(   \bigcup_{j=1}^{\infty}
    B_Y\big(f_j(p_1),\delta'\big)            \Big)       \Big)\leq \mu_Y
    \Big(  \bigcup_{j=1}^{\infty}  \partial B_Y\big(f_j(p_1),\delta'\big)               \Big)=0
    \end{align*}and $\{ (f_n)_{\ast}(\mu_X)    \}_{n=1}^{\infty}$
   converges weakly to $\mu_Y$, we have
   \begin{align*}
    \lim_{n\to \infty}\sum_{j=1}^{\infty} \mu_X \Big( f_n^{-1}\big(
     B_Y(f_j(p_1),\delta')\big)      \Big)= \sum_{j=1}^{\infty}\mu_Y
    \big(  B_Y \big(   f_j(p_1), \delta'             \big)    \big)
    \end{align*}and 
   \begin{align*}
    \lim_{n\to \infty} \mu_X  \Big( f_n^{-1}\big(
     B_Y(f_j(p_1),\delta')\big)      \Big) =   \mu_Y
    \big(  B_Y \big(   f_j(p_1), \delta'             \big)    \big) 
    \end{align*}for any $j$. For any $\varepsilon>0$, there exists
   $k_0\in \mathbb{N}$ such that 
   \begin{align*}
    \sum_{j=1}^{k_0} \mu_Y \big(  B_Y(f_j(p_1), \delta')
    \big) +\varepsilon > \sum_{j=1}^{\infty} \mu_Y \big(  B_Y(f_j(p_1), \delta')         \big).
    \end{align*}Take $n_0 \in \mathbb{N}$ such that 
   \begin{align*}
           \Big|  \sum_{j=1}^{k_0}   \mu_Y \big( B_Y(f_j(p_1 ), \delta')
    \big) - \sum_{j=1}^{k_0}   \mu_X \Big( f_n^{-1}\big(B_Y(f_j(p_1 ), \delta')
    \big)\Big) \Big|< \varepsilon
    \end{align*}and
   \begin{align*}
    \Big|  \sum_{j=1}^{\infty}  \mu_Y \big( B_Y(f_j(p_1 ), \delta')
    \big) -     \sum_{j=1}^{\infty}   \mu_X \Big( f_n^{-1}\big(B_Y(f_j(p_1 ), \delta')
    \big)\Big)                                  \Big| < \varepsilon
    \end{align*}for any $n\geq n_0$. Hence, for any $n\geq n_0$ we
   have
   \begin{align*}
    \sum_{j=1}^{k_0} \mu_X \Big( f_n^{-1} \big(   B_Y (f_j(p_1),
    \delta')         \big)            \Big) + 3\varepsilon >
    \sum_{j=1}^{\infty}\mu_X \Big( f_n^{-1} \big(   B_Y (f_j(p_1),
    \delta')         \big)            \Big),
   \end{align*}which implies that 
   \begin{align*}
    \mu_X \Big( f_n^{-1} \big(   B_Y (f_n(p_1),
    \delta')         \big)            \Big) \to 0 \ \text{as } n\to \infty.
    \end{align*} Fix
   $\delta''>0$ with $\delta'' <\delta'$. Since $p_1 \in X_0$, we get $\dist_Y \big( f_n(p_1),f_n
   (q)\big)\leq \lambda \dist_X(p_1,q) + \varepsilon_n$ for any $q\in
   X_n$ and for any suffieciently large $n\in \mathbb{N}$. Therefore, we
   get
   \begin{align*}
    B_X\Big(     p_1, \frac{\delta''}{\lambda}           \Big) \cap X_n
    \subseteq f_n^{-1}\big( B_Y (f_n (p_1), \delta')\big)
    \end{align*}for any suffieciently large $n\in \mathbb{N}$. Hence, we
   obtain 
   \begin{align*}
    \mu_X   \Big(  B_X   \Big(p_1, \frac{\delta''}{\lambda}    \Big)
    \cap X_n          \Big) \to 0 \ \text{as } n \to \infty,
    \end{align*}which yields $\mu_X \big(B_X  (p_1, \delta''/ \lambda)
   \big) = 0   $. This is a contradition, since $p_1 \in X= \supp
   \mu$. This completes the proof of the claim.
  \end{proof}
  \end{claim}
By virtue of Claim \ref{mikihayoi} and the diagonal argument, we have
 that $\{ f_n(p_j)    \}_{n=1}^{\infty}$ is convergent sequence in $Y$
 for each $j\in \mathbb{N}$. We put $f(p_j):= \lim\limits_{n\to \infty}
 f_n(p_j)$ for any $j\in \mathbb{N}$. Extend the map $f:\{p_1,p_2, \cdots
 \}\to Y$ to $\widetilde{f}:X_0 \to Y$, by using $f$ is a
 $\lambda$-Lipschitz map.
 \begin{claim}\label{harisennbonnnikyouhu}
  For any $\varepsilon >0$, we have $\mu_X \big(\big\{  x\in X \mid \dist_Y
  \big(f_n(x), \widetilde{f}(x)\big)    \geq \varepsilon
  \big\}\big)\to 0$ as $n\to \infty$.
  \begin{proof}Since $X_0 \subseteq
   \bigcup\limits_{j=1}^{\infty}B_X(p_j, \varepsilon /2)$, for any
   $\delta >0$ there exists $k_0\in \mathbb{N}$ such that 
   \begin{align*}
    \mu_X \Big(  \bigcup\limits_{j=1}^{k_0}B_X(p_j, \varepsilon /2) \cap X_0
    \Big)\geq 1-\delta.        
    \end{align*}From the definition, there exists $n_0\in \mathbb{N}$
   such that $\dist_Y \big(       f_n (p_j), \widetilde{f}  (p_j)
   \big) \leq \varepsilon /3   $ for any $n\geq n_0$ and $j=1,2, \cdots,
   k_0$. Take any $x \in \bigcup\limits_{j=1}^{k_0}B_X(p_j, \varepsilon
   /2) \cap X_0$. There exists $1\leq j\leq k_0$ such that $x\in
   B_X(p_j, \varepsilon /2)$. Hence, for any $n\geq n_0$ we have
   \begin{align*}
    \dist_Y \big(  f_n(x) , \widetilde{f}(x)      \big)\leq \ &\dist_Y
    \big( f_n(x),f_n(p_j)     \big) + \dist_Y(f_n(p_j),
    \widetilde{f}(p_j))+ \dist_Y \big(  \widetilde{f}  (p_j),
    \widetilde{f}(x)  \big)\\
    < \ & (\lambda \dist_X(x,p_j)+ \varepsilon_n) +
    \frac{\varepsilon}{3} + \lambda \frac{\varepsilon}{2}\\
    \leq \ &  \varepsilon_n+\frac{\varepsilon}{3}+\lambda \varepsilon. 
    \end{align*}Therefore, for any suffieciently large $n\in
 \mathbb{N}$, we obtain
 \begin{align*}
  \mu_X          \Big( \Big\{  x\in X \mid \dist_Y     \big(f_n(x),
  \widetilde{f}(x)\big)      > \lambda \varepsilon
  +\frac{\varepsilon}{2}        \Big\}            \Big) \leq \mu_X \Big(
  X \setminus \bigcup_{j=1}^{k_0}B_X \Big( p_j, \frac{\varepsilon}{2}
  \Big) \cap X_0 \Big) \leq \delta.
  \end{align*}This completes the proof of the claim.
   \end{proof}
  \end{claim}According to Claim \ref{harisennbonnnikyouhu}, we have
 $\me_1 (f_n , \widetilde{f})\to 0$ as $n\to \infty$. This completes the
 proof of the proposition. 
\end{proof} 
\end{prop}

Gromov proved in \cite[Section $3.\frac{1}{2}.10$]{gromov} the following
proposition by using the distance function $\tra_{\lambda}$ on the space
of finite Borel measures. Although the distance function
$\tra_{\lambda}$ does not appare in the proof of the following
proposition, the proof is essentially the same spirit of his proof.  
\begin{prop}[{cf.~\cite[Section $3.\frac{1}{2}.10$]{gromov}}]\label{saitensinakereba}Let $\{ \mu_n   \}_{n=1}^{\infty}$ be a sequence of Borel
 measures on a metric space $X$ and assume that $\{  \mu_n
 \}_{n=1}^{\infty}$ converges weakly to a Borel measure $\mu$. Then, we
 have 
\begin{align*}
\sikakuu \big(  (X, \dist_X, \mu_n), (X, \dist_X, \mu)       \big)\to 0
 \ \text{as } n\to \infty.
\end{align*}
 \begin{proof}Without loss of generality, we may assume that $\mu(X)=1$
  and $\mu_n(X)=1$ for any $n\in \mathbb{N}$. For any $\varepsilon >0$,
  there exists a sequence $\{  A_i    \}_{i=1}^{\infty}$ of pairwise
  disjoint Borel
  subsets of $X$ satisfying the following properties $(1)-(3)$.
  \begin{itemize}
 \item[$(1)$]$X= \bigcup\limits_{i=1}^{\infty} A_i =X$.
 \item[$(2)$]For any $i\in \mathbb{N}$, $\diam A_i \leq \varepsilon$.
 \item[$(3)$]For any $i,n\in \mathbb{N}$, $\mu(\partial
	     A_i)=\mu_n(\partial A_i)=0$.
   \end{itemize}
From $(1)$ and $(3)$, there exists $m\in \mathbb{N}$ such that
  $\mu\big(\bigcup\limits_{i=1}^{m}A_i\big)= \mu
 \big(\bigcup\limits_{i=1}^{m}\bar{A_i}\big)>1-\varepsilon$. From the
  assumption, $\mu_n(\bar{A_i})=\mu_n(A_i)\to
 \mu(A_i)=\mu(\bar{A_i})$ as $n\to \infty$ for any $i\in
  \mathbb{N}$. Hence, putting 
\begin{align*}
&I_{1n}:= [0,\mu_n(\bar{A_1})), \\
&I_{in}:= 
\Big[\sum_{k=1}^{i-1}\mu_n(\bar{A_k}),\sum_{k=1}^{i}\mu_n(\bar{A_k})\Big),
 \ 
 i=2,3, \cdots, \\
&I_1:= [0,\mu(\bar{A_1})), \\
&I_i:= \Big[\sum_{k=1}^{i-1}\mu(\bar{A_k}),\sum_{k=1}^{i}\mu(\bar{A_k})\Big),\
 i=2,3,\cdots, 
\end{align*}
there exists $N\in \mathbb{N}$ such that
\begin{align*}
\mathcal{L}(I_{in}\cap
 I_i)\geq \mu(\bar{A_i})-\varepsilon /m
\end{align*}for any $n\geq N$ and $i=1,2, \cdots,m$. Fix a parameter $\phi_i:I_i \to
  \bar{A_i}$ of the mm-space $(\bar{A_i},\dist_X,\mu)$ for each $i=1,2, \cdots,m$.
 For any $A\subseteq X$, we indicate by $\nai A$ its interior. Since
  $\mu(\bar{A_i})=\mu(\nai \bar{A_i})$, we have
\begin{align*}
 \mu \Big(\bigcup_{i=1}^m \nai \bar{A_i}\Big)=\sum_{i=1}^m \mu (\bar{A_i})=
 \sum_{i=1}^m \mathcal{L}(I_i)=\mathcal{L}\Big(\bigcup_{i=1}^m I_i\Big).
\end{align*}Take a paramter $\phi :[0,1]\setminus \bigcup\limits_{i=1}^m
  I_i \to X\setminus \bigcup\limits_{i=1}^m \nai \bar{A_i}$ of the
  mm-space $\big(X\setminus \bigcup\limits_{i=1}^m \nai \bar{A_i},
 \dist_X,\mu\big)$. Defining a Borel measurable map $\varphi :[0,1]\to
  X$ by 
\begin{align*}
\varphi(t):=
\left\{
\begin{array}{ll}
\phi_{i}(t) &  \ \ t\in I_i, \ i=1,2, \cdots ,m,\\
\phi(t) &  \ \ t\in [0,1]\setminus
 \bigcup\limits_{i=1}^m I_i, \\
\end{array}
\right.
\end{align*}we see that the map $\varphi$ is a parameter of
  $(X,\dist_X,\mu)$. We take any $n\geq N$. Take parameters 
  $\psi_{in}:I_{in}\to
  \bar{A_i}$ of $i=1,2, \cdots,m$, of the mm-spaces
  $(\bar{A_i},\dist_X,\mu_n)$, and a parameter $\psi_n :[0,1]\setminus
 \bigcup\limits_{i=1}^m I_{in} \to X\setminus \bigcup\limits_{i=1}^m
 \nai \bar{A_i}$ of the mm-space $\big(X\setminus \bigcup\limits_{i=1}^m \nai
  \bar{A_i},\dist_X,\mu_n\big)$. We define a Borel measurable map $\varphi_n :[0,1]\to
 X$ by 
\begin{align*}
\varphi_n(t):=
\left\{
\begin{array}{ll}
\psi_{in}(t) &  \ \ t\in I_{in}, \ i=1,2,\cdots ,m,\\
\psi_n(t) &  \ \ t\in [0,1]\setminus
 \bigcup\limits_{i=1}^m I_{im}. \\
\end{array}
\right.
\end{align*}
The map $\varphi_n$ is a parameter of the mm-space $(X,\dist_X, \mu_n)$
  for each $n\geq N$. Putting $B_n :=\bigcup\limits_{i=1}^{m}(I_i \cap
 I_{in})$, we have
\begin{align*}
\mathcal{L}(B_n)=\sum_{i=1}^m \mathcal{L}(I_i\cap I_{in})
\geq \ &\sum_{i=1}^m(\mu(A_i)-\varepsilon /m)
=  \sum_{i=1}^m \mu(A_i)-\varepsilon \\
= \ & \mu\Big(\bigcup_{i=1}^mA_i\Big)-\varepsilon \geq 1-2\varepsilon.
\end{align*}For any $s,t\in  B_n$, there exist $j,k\in \mathbb{N}$ such
  that $1\leq j,k\leq m$, $s\in I_j\cap I_{jn}$, and $t\in I_k \cap
  I_{kn}$. Since $\varphi(s),\varphi_n(s)\in \bar{A_j}$, $\varphi(t),\varphi_n(t)\in
 \bar{A_k}$, and $(2)$, we have 
\begin{align*}
 \big|\dist_X\big(\varphi(s),\varphi(t)\big)-\dist_X\big(\varphi_n(s),\varphi_n(t)\big)\big|\leq
 \dist_X\big(\varphi(s),\varphi_n(s)\big)+\dist_X
 \big(\varphi(t),\varphi_n(t)\big)\leq 2\varepsilon.
\end{align*}Therefore, we obtain
 $\underline{\square}_1\big((X,\dist_X,\mu_n),(X,\dist_X,\mu)\big)\leq
 \square_1(\varphi_n^{\ast}\dist_X,\varphi^{\ast}\dist_X)\leq
  2\varepsilon$. This completes the proof.
  \end{proof}
\end{prop}

\begin{thm}[{Gromov, cf.~\cite[Section $3\frac{1}{2}.15$, $(3_b')$]{gromov}}]\label{totemomendou}$\sikakuu (X_n, X)\to 0$ as $n\to \infty$ if and only if
 for any $n\in \mathbb{N}$ there exist a Borel measurable map $ p_n:X_n
 \to X  $, a Borel subset $\widetilde{X}_n \subseteq X_n$, and a
 positive number $
 \varepsilon_n $ satisfying the
 following conditions $(1)-(4)$.
\begin{itemize}
                   \item[$(1)$]$\varepsilon_n \to 0$ as $n\to \infty$.
		   \item[$(2)$]$\mu_{X_n}(X_n \setminus
			       \widetilde{X}_n)\leq \varepsilon_n$ for
			       $n=1,2, \cdots$.
		   \item[$(3)$]$|\dist_{X_n}(x,y)-\dist_X\big(p_n(x),p_n(y)\big)|\leq
			       \varepsilon_n$ for any $x,y\in
			       \widetilde{X}_n$.
		   \item[$(4)$]The sequence $\{  (p_n)_{\ast}(\mu_{X_n})
			       \}_{n=1}^{\infty}$ converges weakly to
			       $\mu_X$.
\end{itemize}
 \begin{proof}
Assume that $(1)-(4)$ holds. By virtue of Proposition
  \ref{saitensinakereba}, we have $\sikakuu (X_n, X)\to 0$ as $n\to \infty$.

Assume that $\sikakuu (X_n,X)\to 0$ as $n\to \infty$. Without loss of
  generality, we may assume that $\mu_X(X)=\mu_{X_n}(X_n)=1$ for any
  $n\in \mathbb{N}$. From the assumption, there exist parameters
  $\varphi:[0,1]\to X$ of $X$ and $\varphi_n:[0,1]\to X_n$ of $X_n$,
  $n\in \mathbb{N}$, such that $\square_1(\varphi_n^{\ast}\dist_{X_n},
  \varphi^{\ast}\dist_X) \to 0$ as $n \to \infty$. Hence, for each
  $n=1,2, \cdots$, there exist $\varepsilon_n>0$ and compact subset $K_n\subseteq
	  [0,1]$ satisfying the following conditions $(1)'-(4)'$: 
	  \begin{itemize}
	   \item[$(1)'$]$\varepsilon_n \to 0$ as $n\to \infty$.
	   \item[$(2)'$]$\mathcal{L}(K_n)> 1-\varepsilon_n$.
	   \item[$(3)'$]For any $s,t \in K_n$, 
	    $\big|\dist_X\big(\varphi(s),\varphi(t)\big)-\dist_{X_n}\big(\varphi_n(s),\varphi_n(t)\big)\big|<\varepsilon_n$. 
	   \item[$(4)'$]The maps $\varphi|_{K_n}:K_n\to
	    X$ and $\varphi_n|_{K_n}:K_n \to X_n$ are continuous. 
	  \end{itemize}
	  By $(4)'$, each set $\varphi_n(K_n)$ is compact. For each
  $n\in \mathbb{N}$, there exist $l_n\in \mathbb{N}$ and a sequence $\{ B_{in}
  \}_{i=1}^{l_n}$ of pairwise disijoint Borel subsets of $X_n$ such that
  $\diam B_{ni} < \varepsilon_n$ for any $i$ and 
	  $\varphi_n(K_n)=\bigcup\limits_{i=1}^{l_n}B_{in}$. For each
  $i$, we fix a point $p_{in}\in B_{in}$. Then there exist a point
  $t_{in}\in K_n$ with $p_{in}=\varphi_n(t_{in})$. Put
  $q_{in}:=\varphi(t_{in})\in X$.
	  \begin{claim}$\varphi(K_n)\subseteq
	   \bigcup\limits_{i=1}^{l_n}B_X(q_{in},2\varepsilon_n)$. 
	   \begin{proof}Take any $q= \varphi(s)\in
	  \varphi(K_n)$ with $s\in K_n$. Since $\varphi_n(s)\in
	    \varphi_n(K_n)\subseteq \bigcup\limits_{i=1}^l
	    B_{X_n}(p_{in},\varepsilon_n)$, there exists $1\leq i \leq
	  l_n$ such that
	    $\dist_{X_n}\big(\varphi_n(s),\varphi_n(t_{in})\big)<
	    \varepsilon_n$. Hence, by $(3)'$, we have  
	  \begin{align*}
	  \dist_X(q,q_{in})=\dist_X(\varphi(s),\varphi(t_{in}))< \dist_{X_n}(\varphi_n(s),\varphi_n(t_{in}))+\varepsilon_n <2\varepsilon_n.
	  \end{align*}This completes the proof of the claim.
	    \end{proof}
	   \end{claim}We denote by
  $\widetilde{q}_{1n},\widetilde{q}_{2n},\cdots ,\widetilde{q}_{m_n n}$
  the mutually different elements of $\{q_{1n},q_{2n},\cdots ,q_{l_n n}\}$. Put
	  \begin{align*}
	   &C_{1n}:=\varphi(K_n)\cap B_X(\widetilde{q}_{1n},2\varepsilon_n
	   )\setminus \{ \widetilde{q}_{2n},\widetilde{q}_{3n},\cdots ,\widetilde{q}_{m_n n}    \}, \\
	   &C_{in}:= \varphi(K_n)\cap
	   B_X(\widetilde{q}_{in},2\varepsilon_n) \setminus
	   \Big\{ \bigcup_{j=1}^{i-1}\big(B_X(\widetilde{q}_{jn},2\varepsilon_n)\setminus \{ \widetilde{q}_{in}
	   \} \big)\cup
	   \{\widetilde{q}_{i+1 n},\widetilde{q}_{i+2n},\cdots , \widetilde{q}_{m_n n}\} \Big\},\\
	   &\hspace{12.2cm} i=2,3,\cdots , m_n.
	  \end{align*}
	  It is easy to see that $\widetilde{q}_{in}\in C_{in},\ \varphi(K_n)=
	  \bigcup\limits_{j=1}^{m_{n}}C_{jn}$, $C_{in}\cap
	  C_{jn}=\emptyset$ for $i\neq j$, and $\diam C_{in}\leq
  4\varepsilon_n$. Take points $x_n^0 \in X_n$ for any $n\in
  \mathbb{N}$ and $x^0\in X$. We define a Borel measurable map 
	  $p_{n}:X_n \to X$ by $p_n(x_n):=q_{in}$ if $x_n\in B_{in}$ and
  $p_n(x_n):=x^0$ if $x_n\in X_n \setminus \varphi_{n}(K_n)$. For each
  $i=1,2, \cdots ,m_n$, we fix $j$ with $\widetilde{q}_{in}=q_{jn}$ and
  put $k_n(i):=j$. 

	  \begin{claim}
	   The sequence $\{ (p_n)_{ \ast}(\mu_{X_n} ) \}_{n=1}^{\infty}$
	   converges weaky to the measure $\mu_X$.
	   \begin{proof}Let $g:X\to \mathbb{R}$ be any bounded
	    uniformly continuous function and put $M:=\sup\limits_{x\in
	    X}|g(x)|$. We shall prove
	  \begin{align*}
	  \int_{X_n}(g\circ p_n)(x_n) \ d\mu_{X_n}(x_n)\ \to \ \int_X g(x)\ d\mu_X(x)\
	  \ \text{as }n\to \infty.
	  \end{align*}Since 
	  \begin{align*}
	   \int_{X_n}(g\circ p_n)(x_n)\ d\mu_{X_n}(x_n)
	   =\ &\int_{0}^{1}(g\circ p_n \circ
	   \varphi_n)(s)\ d\mathcal{L}(s)\\
	   =\ &\int_{K_n}(g\circ p_n \circ \varphi_n)(s)\
	   d\mathcal{L}(s)+\int_{[0,1]\setminus K_n}(g\circ p_n \circ \varphi_n)(s)\ d\mathcal{L}(s),
	  \end{align*}we get
	  \begin{align*}
	   &\Big| \int_{X_n}(g\circ p_n)(x_n)\ d\mu_{X_n}(x_n)-\int_{K_n}(g\circ p_n\circ \varphi_n)(s)\ d\mathcal{L}(s)\Big|\\ \leq &\
	   \int_{[0,1]\setminus K_n}|(g\circ p_n \circ \varphi_n)(s)|\ d\mathcal{L}(s)
	   \hspace{0.2cm} \leq \  M\varepsilon_n.
	  \end{align*}Similary, we have 
	  \begin{align*}
	   \Big| \int_X g(x)\ d\mu_X(x) -\int_{K_n}(g\circ
	   \varphi)(s)\ d\mathcal{L}(s)          \Big|\leq M\varepsilon_n.
	  \end{align*}Since for any $s\in \varphi_n^{-1}(B_{in})\cap
	  \varphi^{-1}(C_{jn})$ 
	  \begin{align*}
	   \dist_{X_n}(\varphi_n(s),p_{in})\leq
	  \varepsilon_n \text{ and }\dist_X(\varphi(s),\widetilde{q}_{jn})\leq 2\varepsilon_n,
	  \end{align*}we obtain 
	  \begin{align*}
	   &\dist_X(\varphi(s),q_{in})\\ \leq \ &
	   \dist_X(\varphi(s),\widetilde{q}_{jn})+|\dist_X(\widetilde{q}_{jn},q_{in})-\dist_{X_n}(\varphi_n(s),p_{in})|+\dist_{X_n}(\varphi_n(s),p_{in})\\
	   <\ & 2\varepsilon_n +
	   |\dist_X(\widetilde{q}_{jn},q_{in})-\dist_{X_n}(p_{k_n(j)n},p_{in})|+\dist_{X_n}(p_{k_n(j)n},\varphi_n(s))
	   +\varepsilon_n \\
	   <\ & 4\varepsilon_n + \dist_{X_n}(p_{k_n(j)n},\varphi_n(s))\\
	   <\ & 4\varepsilon_n + \dist_X(\widetilde{q}_{jn},\varphi(s))+\varepsilon_n \ \leq 7\varepsilon_n. 
	  \end{align*}Since $g$ is uniformly continuous function on $X$, for any
	    $\varepsilon > 0$ there exists $\delta>0$ such that
	    $|g(x)-g(y)|< \varepsilon$ for any
	    $x,y\in X$ with $\dist_X(x,y)<\delta$. Hence for any $n\in
	    \mathbb{N}$ with $7\varepsilon_n< \delta$, we have
	    $|g(q_{in})-g(\varphi(s))|< \varepsilon$, which implies that 
	  \begin{align*}
	   &\Big| \int_{K_n}(g\circ p_n \circ \varphi_n)(s)\
	   d\mathcal{L}(s)-\int_{K_n}(g\circ \varphi)(s)\
	   d\mathcal{L}(s)
	   \Big|\\
	   \leq \ & \sum_{i,j=1}^{l_n,m_n}\int_{\varphi_{n}^{-1}(B_{in})\cap
	   \varphi^{-1}(C_{jn})\cap K_n} |g(q_{in})-g(\varphi (s))|\ d\mathcal{L}(s)\\
	   <\ & \varepsilon \sum_{i,j=1}^{l_n,m_n}
	   \mathcal{L}(\varphi_n^{-1}(B_{in})\cap
	   \varphi^{-1}(C_{jn})\cap K_n)=\varepsilon\mathcal{L}(K_n)\leq
	   \varepsilon. 
	  \end{align*}
	  Therefore, we obtain
	  \begin{align*}
	   &\Big| \int_{X_n} (g\circ p_n)(x_n) \ d\mu_{X_n}(x_n) -\int_{X}g(x)\
	   d\mu_{X}(x) \Big|\\
	   \leq \ & \Big|  \int_{X_n}(g\circ p_n)(x_n) \ d\mu_{X_n}(x_n) -
	   \int_{K_n} (g\circ p_n \circ \varphi_n)(s)\
	   d\mathcal{L}(s) \Big|\\
	   & +\Big| \int_{K_n}(g\circ p_n\circ \varphi_n)(s) \
	   d\mathcal{L}(s)   -  \int_{K_n} (g\circ \varphi)(s)\
	   d\mathcal{L}(s)\Big|\\
	   & +\Big|  \int_{K_n}(g\circ \varphi)(s)\ d\mathcal{L}(s)-
	   \int_{X}g(x)\ d\mu_X(x)   \Big|\\
	   \leq \ & 2M\varepsilon_n +\varepsilon.
	  \end{align*}This completes the proof of the claim.
	   \end{proof}
	   \end{claim}

	  For any $x\in B_{in},y\in B_{jn}$, we obtain
	   \begin{align*}
		&|\dist_{X_n}(x,y)-\dist_X(p_n(x),p_n(y))|\\=
		\ &|\dist_{X_n}(x,y)-\dist_X(q_{in},q_{jn})|\\
		\leq \ & |\dist_{X_n}(x,y)-\dist_{X_n}(p_{in},p_{jn})|+|\dist_{X_n}(p_{in},p_{jn})-\dist_X(q_{in},q_{jn})|\\
		\leq \ & \dist_{X_n}(x,p_{in})+\dist_{X_n}(y,p_{jn})+|\dist_{X_n}(p_{in},p_{jn})-\dist_X(q_{in},q_{jn})|\\
		< \  & 2\varepsilon_n+ \varepsilon_n\\
		=\ & 3\varepsilon_n.
	   \end{align*}Therefore, we have complete the proof of Theorem \ref{totemomendou}. 
  \end{proof}
\end{thm}

Modifying the proof of Theorem \ref{totemomendou}, we get the following
corollary:
\begin{cor}Let $X$ and $X_n$, $n\in \mathbb{N}$, be compact mm-spaces. Assume that
$X=\supp \mu_X$, $X_n = \supp \mu_{X_n}$, and
 $\mu_X(X)= \mu_{X_n}(X_n)$ for any $n\in \mathbb{N}$.
Then, the sequence $\{X_n\}_{n=1}^{\infty}$ converges to $X$ with
 respect to $\underline{\square}_0$ if and only if $\{ X_n
 \}_{n=1}^{\infty}$ converges to $X$ in the sense of the measured
 Gromov-Hausdorff convergence.
\end{cor}

Combining Proposition \ref{kumogaippai} and Theorem \ref{totemomendou}, we get
the following corollary:
\begin{cor}Assume that $\sikaku (X,Y)=0$. Then, two mm-spaces $X$ and
 $Y$ are isomorphic to each other. 
\end{cor}

\section{Stability of homogenuity}
We say that an mm-space $X$ \emph{Lipschitz dominates} an mm-space $Y$
and write $X\succ Y$ if
there exist $1$-Lipschitz map $p:\supp \mu_X \to \supp \mu_Y$ and $c\geq
1$ such that $p_{\ast}(\mu_X)= c\mu_Y$.

\begin{thm}[{Gromov, cf.~\cite[Section $3\frac{1}{2}.15$, $(b)$]{gromov}}]\label{amefurisou}Assume that $\sikaku (X_n, X)$, $\sikaku(Y_n,Y)\to 0$ as
 $n\to \infty$ and $X_n \succ Y_n$ for any $n\in \mathbb{N}$. Then we
 have $X\succ Y$. 
\begin{proof}Without loss of generality, we may assume that
 $\mu_{X_n}(X_n)=\mu_{Y_n}(Y_n)= \mu_{X}(X)=\mu_Y(Y)=1$, $X=\supp
 \mu_X$, $Y=\supp \mu_Y$, $X_n= \supp \mu_{X_n}$, and $Y_n=\supp \mu_{Y_n}$ for any $n\in
 \mathbb{N}$. From the assumption, for any $n\in \mathbb{N}$ there exists
 a $1$-Lipschitz map $f_n:X_n \to Y_n$ such that
 $(f_n)_{\ast}(\mu_{X_n})=\mu_{Y_n}$. By using Theorem \ref{totemomendou}, for
 any $n\in \mathbb{N}$
 there exists a Borel measurable map $q_n:Y_n \to Y$, a compact subset
 $\widetilde{Y}_n \subseteq Y_n$, and $\varepsilon_n >0$ such that
\begin{itemize}
                   \item[$(1)$]$\varepsilon_n \to 0$ as $n\to \infty$,
		   \item[$(2)$]$\mu_{Y_n}(Y_n \setminus
			       \widetilde{Y}_n)\leq \varepsilon_n$ for
			       $n=1,2, \cdots$,
		   \item[$(3)$]$|\dist_{Y_n}(x,y)-\dist_Y\big(q_n(x),q_n(y)\big)|\leq
			       \varepsilon_n$ for any $x,y\in
			       \widetilde{Y}_n$,
		   \item[$(4)$]The sequence $\{  (q_n)_{\ast}(\mu_n)
			       \}_{n=1}^{\infty}$ converges weakly to
			       $\mu_Y$.
\end{itemize}From now on, we define a Borel measurable map $p_n :X\to
 X_n$ as follows: Since $\sikakuu(X_n,X)\to 0$ as $n\to \infty$, there
 exists a parameter $\varphi_n:[0,1]\to X_n$ and $\varphi
 :[0,1]\to X$ such that $\square_1 (\varphi_n^{\ast}\dist_{X_n},
 \varphi^{\ast} \dist_X)\to 0$ as $n\to \infty$. Hence, there exists a
 compact subset $K_n \subseteq [0,1]$ and $\varepsilon_n' >0$ satisfying
 the following properties $(1)'-(4)'$:
 \begin{itemize}
  \item[$(1)'$]$\varepsilon_n'\to 0$ as $n\to \infty$.
  \item[$(2)'$]$\mathcal{L}(K_n)>1-\varepsilon_n'$.
  \item[$(3)'$]For any $s,t \in K_n$,
	       $\big|\dist_{X_n}\big(\varphi_n(s),
	       \varphi_n(t)\big)-\dist_X \big(\varphi(s),
	       \varphi(t)\big)\big|< \varepsilon_n'$.
  \item[$(4)'$]The maps $\varphi_n|_{K_n}:K_n \to [0,1]$ and
	       $\varphi_n|_{K_n}:K_n \to [0,1]$ are continuous.
  \end{itemize}By $(4)'$, the sets $\varphi_n(K_n)\cap
 f_n^{-1}(\widetilde{Y}_n)$ and $\widetilde{X}_n:= \varphi \big(K_n \cap
 \varphi_n^{-1}\big(f_n^{-1}(\widetilde{Y}_n)\big)\big)$ are compact. For each
  $n\in \mathbb{N}$, there exist $l_n\in \mathbb{N}$ and a sequence $\{ C_{in}
  \}_{i=1}^{l_n}$ of pairwise disijoint Borel subsets of $X$ such that
  $\diam C_{in} < \varepsilon_n$ for any $i$ and 
	  $\widetilde{X}_n=\bigcup\limits_{i=1}^{l_n}C_{in}$. For each
  $i$, we fix a point $q_{in}\in C_{in}$. Then there exist a point
  $t_{in}\in K_n \cap \varphi_n^{-1}\big( f_n^{-1}(\widetilde{Y}_n)
 \big)$ with $q_{in}=\varphi(t_{in})$. Put $p_{in}:=\varphi(t_{in})\in
 X$. Then, we get $\varphi_n (K_n)\cap f_n^{-1}(\widetilde{Y}_n)
 \subseteq \bigcup\limits_{i=1}^{l_n}B_{X_n}(p_{in},2\varepsilon_n')$. We denote by
  $\widetilde{p}_{1n},\widetilde{p}_{2n},\cdots ,\widetilde{p}_{m_n n}$
  the mutually different elements of $\{p_{1n},p_{2n},\cdots ,p_{l_n n}\}$. Put
	  \begin{align*}
	   &B_{1n}:=\varphi_n(K_n) \cap f_n^{-1}(\widetilde{Y}_n)\cap B_{X_n}(\widetilde{p}_{1n},2\varepsilon_n'
	   )\setminus \{ \widetilde{p}_{2n},\widetilde{p}_{3n},\cdots ,\widetilde{p}_{m_n n}    \}, \\
	   &B_{in}:= \varphi_n(K_n)\cap f_n^{-1}(\widetilde{Y}_n) \cap 
	   B_{X_n}(\widetilde{p}_{in},2\varepsilon_n') \\ 
	   & \hspace{2cm} \setminus
	   \Big\{
	   \bigcup_{j=1}^{i-1}\big(B_{X_n}(\widetilde{p}_{jn},2\varepsilon_n')      \setminus \{ \widetilde{p}_{in}
	   \} \big)\cup
	   \{\widetilde{p}_{i+1 n},\widetilde{p}_{i+2n},\cdots ,
	   \widetilde{p}_{m_n n}\} \Big\}, i=2,3,\cdots , m_n.
	  \end{align*}
	  It is easy to see that $\widetilde{p}_{in}\in B_{in},\
 \varphi_n(K_n)\cap f_n^{-1}(\widetilde{Y}_n)=
	  \bigcup\limits_{j=1}^{m_{n}}B_{jn}$, $B_{in}\cap
	  B_{jn}=\emptyset$ for $i\neq j$, and $\diam B_{in}\leq
  4\varepsilon_n'$. Take points $x_n^0 \in X_n$ for any $n\in
  \mathbb{N}$ and $x^0\in X$. We put $p_n(x):=p_{in}$ if $x\in C_{in}$ and
  $p_n(x):=x^0$ if $x\in X \setminus \widetilde{X}_n$. 
The same proof in Theorem \ref{totemomendou} implies the following: There
 exists a positive number $\delta_n >0$ such that $\delta_n \to 0$ as
 $n\to \infty$, $\mu_X(X\setminus
 \widetilde{X}_n)< \delta_n$, and $\big|\dist_{X_n}\big(  p_n(x),p_n(x')
 \big)- \dist_X(x,x')\big|< \delta_n$ for any $x,x'\in \widetilde{X}_n$.

Put $g_n:=q_n \circ f_n \circ p_n:X\to Y$. For any $x,x' \in
 \widetilde{X}_n$,  
\begin{align*}
\dist_Y\big(g_n(x),g_n(x')\big)-\dist_X(x,x')\leq \ &\big|\dist_Y
 \big(g_n(x),g_n(x')\big)- \dist_{Y_n} \big(  (f_n\circ p_n)(x), (f_n
 \circ p_n)(x')     \big)            \big|\\
 \ & \hspace{1.5cm}+ \dist_{Y_n} \big(  (f_n\circ p_n)(x), (f_n
 \circ p_n)(x')     \big)   -\dist_X(x,x')\\
\leq \ &\varepsilon_n + \big( \dist_{X_n}\big( p_n(x),p_n(x')\big)-\dist_X(x,x') \big)\\
\ &+   \big( \dist_{Y_n} \big(  (f_n\circ p_n)(x), (f_n
 \circ p_n)(x')     \big)-\dist_{X_n}\big(p_n(x),p_n(x') \big)
 \big)\\
\leq \ &\varepsilon_n +\delta_n.
\end{align*}Hence, $g_n$ is a $1$-Lipschitz up to $(\varepsilon_n
 +\delta_n)$ Borel measurable map.
 \begin{claim}\label{atukuteiya}The sequence $\{  (g_n)_{\ast}(\mu_X)
  \}_{n=1}^{\infty}$ converges weakly to the measure $\mu_Y$.
\begin{proof}
Let $h:Y\to \mathbb{R}$ be any bounded uniformly continuous function on
 $Y$. We will prove that
\begin{align*}
\int_X (h \circ g_n )(x) \ d\mu_X(x) \to \int_Y h(y) \ d\mu_Y(y) \
 \text{as }n\to \infty.
\end{align*}Since
 \begin{align*}
\int_{X_n}(h\circ q_n \circ f_n)(x_n)\ d\mu_{X_n}(x_n)=\int_{Y_n}
  (h\circ q_n)(y_n)\ d\mu_{Y_n}(y_n) \to \int_Y h \ d\mu_Y(y) \text{ as
  }n\to \infty,
\end{align*}it suffices to prove that
\begin{align*}
\Big|    \int_X (h\circ g_n)(x)\ d\mu_X(x)- \int_{X_n}(h\circ q_n \circ
 f_n)(x_n) \ d\mu_{X_n}(x_n)                  \Big|\to 0 \text{ as }n\to
\infty.\end{align*}Put $M:= \sup\limits_{y\in Y}|h(y)|$. Take any
 $\varepsilon >0$. We have
\begin{align*}
&\Big| \int_X (h\circ g_n)(x)\ d\mu_X(x)- \int_{K_n \cap
 \varphi_n^{-1}( f_n^{-1}(\widetilde{Y}_n))}   ( h\circ g_n \circ
 \varphi)(s) \ d\mathcal{L}(s)     \Big| \\
\leq \ &\int_{[0,1]\setminus \big(K_n \cap
 \varphi_n^{-1}( f_n^{-1}(\widetilde{Y}_n))\big)} M \ d\mathcal{L}(s)< M\varepsilon 
\end{align*}and 
\begin{align*}
&\Big|    \int_{X_n}(h\circ q_n \circ f_n)(x_n)    \ d\mu_{X_n}(x_n)
 - \int_{K_n \cap \varphi_n^{-1}(f_n^{-1}(\widetilde{Y}_n))} (h \circ
 q_n \circ f_n \circ \varphi_n)(s) \ d\mathcal{L}(s)  \Big|\\
\leq
 \ & \int_{[0,1]\setminus \big(K_n      \cap \varphi_n^{-1}(f_n^{-1}(\widetilde{Y}_n))\big)}M \ d\mathcal{L}(s)< M \varepsilon
\end{align*}for any suffieciently large $n\in \mathbb{N}$. For any $\delta >0$, we put 
\begin{align*}
\rho_h(\delta):= \sup \{   |h(u)-h(v)| \mid \dist_Y(u,v)< \delta ,\ u,v\in
 Y\}.
\end{align*}Let $\varepsilon' >0$ with $\rho_h(2\varepsilon')< \varepsilon$. For any $s\in K_n \cap \varphi_n^{-1}(f_n^{-1}(\widetilde{Y}_n)) \cap
 \varphi^{-1}(C_{kn})\cap \varphi_n^{-1}(B_{jn})$, we get
 $\dist_{X_n}(\varphi_n(s),p_{kn})< \varepsilon'$ for suffieciently large
 $n\in \mathbb{N}$ by the same method of the proof in Theorem
 \ref{totemomendou}. Assume that $x,y \in f_n^{-1}(\widetilde{Y}_n)$ and
 $\dist_{X_n}(x,y)< \varepsilon'$. Then, for any suffieciently large
 $n\in \mathbb{N}$, we have
\begin{align*}
\dist_Y\big((q_n \circ f_n)(x), (q_n \circ f_n)(y)\big)\leq
 \dist_{Y_n}\big(    f_n(x),f_n(y) \big)+ \varepsilon_n \leq
 \dist_X(x,y)+\varepsilon_n <2 \varepsilon'.
\end{align*}Hence, we get 
\begin{align*}
&\Big|   \int_{K_n \cap \varphi_n^{-1}(f_n^{-1}(\widetilde{Y}_n))}\Big( (h\circ
 q_n \circ f_n)\big(  p_n (\varphi(s) )  \big)-   (h \circ q_n \circ
 f_n)\big(    \varphi_n(s)\big)   \Big)\ d\mathcal{L}(s) \Big|\\
\leq \ & \sum_{k,j=1}^{l_n,m_n}\int_{K_n \cap
 \varphi_n^{-1}(f_n^{-1}(\widetilde{Y}_n))  \cap \varphi^{-1}(C_{kn})
 \cap \varphi_n^{-1}(B_{jn})    }\Big|    (h\circ q_n \circ
 f_n)(p_{kn})- (h\circ q_n \circ f_n)\big(\varphi_n(s)\big)   \Big| \
 d\mathcal{L}(s)\\
\leq \ & \sum_{k,j=1}^{l_n,m_n}\int_{K_n \cap
 \varphi_n^{-1}(f_n^{-1}(\widetilde{Y}_n))  \cap \varphi^{-1}(C_{kn})
 \cap \varphi_n^{-1}(B_{jn})   } \rho_h(2\varepsilon')\
 d\mathcal{L}(s) \leq \varepsilon. 
\end{align*}
This completes the proof of the claim.
\end{proof}
  \end{claim}
Combining Proposition \ref{kumogaippai} and Claim \ref{atukuteiya}, we may
 assume that the sequence $\{ g_n  \}_{n=1}^{\infty}$ converges with
 respect to the distance function $\me_1$. Let $g:X\to Y$ be its
 limit. Then this $g$ is obviously a $1$-Lipschitz map.
\begin{claim}\label{gumihaoisii}The sequence $\{ (g_n)_{\ast} (\mu_{X_n})
 \}_{n=1}^{\infty}$ converges weakly to the measure $g_{\ast}(\mu_X)$.
\begin{proof}Let $U\subseteq Y$ be any open subset. Put $U(\delta):=\{
  y\in U \mid \dist_Y(y, X\setminus U)> \delta\}$ for any
 $\delta>0$. For any $\varepsilon >0$, there exists $\delta >0$ such
 that $\mu_X\big(f^{-1}(U)\big)< \mu_X \big(f^{-1}(U(\delta))\big)+
 \varepsilon$. Therefore, we obtain
\begin{align*}
\mu_X \big(f^{-1}(U)\big)<\ &
 \mu_X\big(f^{-1}(U(\delta))\big)\\
=\ &\limsup_{n\to \infty} \mu_X \big(
 f^{-1}(U(\delta)) \cap \{   x \in X \mid \dist_Y (f_n(x),f(x))
< \delta
 \}    \big)\\ \leq \ &\liminf_{n\to \infty} \mu_X\big(f_n^{-1}(U)\big).
\end{align*}This completes the proof of the claim.
\end{proof}
\end{claim}Combining Claim \ref{atukuteiya} and Claim \ref{gumihaoisii},
 we get $g_{\ast}(\mu_X)= \mu_Y$. This completes the proof of the theorem.
\end{proof}
\end{thm}

Modifying the proof of Theorem \ref{amefurisou}, we get the following corollary:
\begin{cor}Assume that a sequence $\{  M_n    \}_{n=1}^{\infty}$ of
 compact homogeneous Riemannian manifolds convergence to an mm-space $X$
 with respect to the distance function $\sikaku$ and $X=\supp \mu_X$. Then, the limit space $X$ is also homogeneous and every
 isometry $g:X\to X$ satisfy $g_{\ast}(\mu_X)=\mu_X$. 
\end{cor}
	
	\end{document}